\documentclass[12pt]{article}
\usepackage{amssymb}
\begin{document}
\baselineskip=20pt
\parskip=5pt
\newcommand{\f}{\frak }
\newcommand{\s}{{\frak s}}
\newcommand{\h}{{\frak h}}
\newcommand{\sll}{{\frak sl}}
\newcommand{\rR}{{\frak r}}
\newcommand{\n}{{\frak n}}
\newcommand{\rr}{{\frak r}^{\bot}}
\newcommand{\cs}{C_{{\frak g}}({\frak s})}
\newcommand{\cg}{C({\frak g})}
\newcommand{\ccr}{C({\frak r})}
\newcommand{\e}{\varepsilon}
\newcommand{\diag}{{\rm diag}}
\newcommand{\Hom}{{\rm Hom\,}}
\newcommand{\ad}{{\rm ad\,}}
\def\dis{\displaystyle}
\newcommand{\rd}{\mbox{Rad}}
\newcommand{\kn}{\mbox{ker}}
\newcommand{\psp}{\vspace{0.4cm}}
\newcommand{\pse}{\vspace{0.2cm}}
\newcommand{\ptl}{\partial}
\newcommand{\dlt}{\delta}
\newcommand{\Dlt}{\Delta}
\newcommand{\sgm}{\sigma}
\newcommand{\al}{\alpha}
\newcommand{\be}{\beta}
\newcommand{\G}{\Gamma}
\newcommand{\gm}{\gamma}
\newcommand{\lmd}{\lambda}
\newcommand{\td}{\tilde}
\newcommand{\stl}{\stackrel}
\newcommand{\ol}{\overline}
\newcommand{\es}{\epsilon}
\newcommand{\la}{\langle}
\newcommand{\ra}{\rangle}
\newcommand{\vf}{\varphi}
\newcommand{\vsi}{\varsigma}
\newcommand{\ves}{\varepsilon}
\newcommand{\vt}{\vartheta}
\newcommand{\wt}{\mbox{wt}\:}
\newcommand{\sym}{\mbox{sym}}
\newcommand{\for}{\mbox{for}}
\newcommand{\mbb}{\mathbb}
\def\qed{\hfill \hfill \ifhmode\unskip\nobreak\fi\ifmmode\ifinner\else\hskip5pt
\fi\fi
 \hbox{\hskip5pt\vrule width4pt height6pt depth1.5pt\hskip 1 pt}}

\def\der{\mbox{der}}
\def\a{\alpha}
\def\b{\beta}
\def\sF{\hbox{$\sc I\hskip -3.5pt F$}}
\def\Z{\hbox{$Z\hskip -5.2pt Z$}}
\def\sZ{\hbox{$\sc Z\hskip -4.2pt Z$}}
\def\Q{\hbox{$Q\hskip -5pt\vrule height 6pt depth 0pt\hskip 6pt$}}
\def\R{\hbox{$I\hskip -3pt R$}}
\def\C{\hbox{$C\hskip -5pt \vrule height 6pt depth 0pt \hskip 6pt$}}
\def\J {\vec J}
\def\d{\delta}
\def\D{\Delta}
\def\g{\gamma}
\def\G{\Gamma}
\def\l{\lambda}
\def\L{\Lambda}
\def\o{\omiga}
\def\p{\psi}
\def\Si{\Sigma}
\def\si{\sigma}
\def\sc{\scriptstyle}
\def\ssc{\scriptscriptstyle}
\def\dis{\displaystyle}
\def\cl{\centerline}
\def\nl{\newline}
\def\DD{{\cal D}}
\def\ll{\leftline}
\def\rl{\rightline}
\def\sF{\hbox{$\sc I\hskip -2.5pt F$}}
\def\ol{\overline}
\def\ul{\underline}
\def\wt{\widetilde}
\def\wh{\widehat}
\def\rar{\rightarrow}
\def\Rar{\Rightarrow}
\def\lar{\leftarrow}
\def\Lar{\Leftarrow}
\def\rla{\leftrightarrow}
\def\Rla{\Leftrightarrow}
\def\bs{\backslash}
\def\hs{\hspace*}
\def\vs{\vspace*}
\def\rb{\raisebox}
\def\ra{\rangle}
\def\la{\langle}
\def\Rad{\mbox{Rad}}
\def\SS{\hbox{$S\hskip -6.2pt S$}}
\def\hi{\hangindent}
\def\ha{\hangafter}
\def\ni{\noindent}
\def\Bbb#1{{\mbox{b$\!\!\!\!$}\cal #1}}
\def\AA{{\cal A}}
\def\BB{{\cal B}}
\def\CC{{\cal C}}
\def\JJ{{\cal J}}
\def\KK{{\cal K}}

\def\v{\vec}
\def\SGN{{\rm sgn}}
\def\HOM{\mbox{Hom}'_{\sZ}(\G,\mbb{F})}
\def\hom{\mbox{Hom}^*_{\sZ}(\G,\mbb{F})}

\def\vi{\vec i}
\def\vj{\vec j}
\def\vk{\vec k}
\def\vm{\vec m}
\def\ii{{\bf i}}
\def\jj{{\bf j}}
\def\kk{{\bf k}}
\def\sone{{1\hskip -5.5pt 1}}
\def\one{{1\hskip -6.5pt 1}}
\def\sJ{J}
\def\th{\theta}
\def\isom{\mbox{${}^{\cong}_{\rar}$}}
\def\N{\mathbb{N}}
\def\Z{\mathbb{Z}}
\def\sZ{\mathbb{Z}}
\def\Q{\mathbb{Q}}
\def\R{\mathbb{R}}
\def\C{\mathbb{C}}
\def\sC{\mathbb{C}}
\def\sF{\mathbb{F}}
\def\F{\mathbb{F}}
\def\fg{{\frak g}}
\def\derg{{\rm Der\,}{\frak g}}
\def\semiprod{\,{\large \times}\hskip-3pt\mbox{\rb{1.5pt}{$\ssc|\,$\,}}}
\par \par \cl{\Large\bf Derivation Algebras of Centerless
Perfect} \cl{{\Large\bf Lie Algebras Are Complete}\footnote{This
work is supported by the NSF of China, two grants ``Excellent Young
Teacher Program'' and ``Trans-Century Training Programme Foundation
for the Talents'' from Ministry of Education of China, the Doctoral
Programme Foundation of Institution of Higher Education and a fund
of Jiangsu Educational Committee.}} \vskip4pt \cl{(appeared in {\it
J. Algebra}, {\bf285} (2005), 508--515.)}
\par\vs{-9pt}\ \par
\cl{Yucai Su$^*$, \ Linsheng Zhu$^\dag$}
\par\vs{-9pt}\ \par
{\small\it $^*$ Department of Mathematics, Shanghai Jiaotong
University, Shanghai 200030, China\par $^\dag$ Department of
Mathematics,
 Changshu Institute of Technology, Jiangsu 215500, China}
\par\vs{-9pt}\ \par
{\bf Abstract:} \ It is proved that the derivation algebra of a
centerless perfect Lie algebra of arbitrary dimension over any
field of arbitrary characteristic is complete and that the
holomorph of a centerless perfect Lie algebra is complete if and
only if its outer derivation algebra is centerless.
\par\vs{-10pt}\ \par
{\bf Key works:} \ Derivation, complete Lie algebra, holomorph of
Lie algebra
\par\vs{-10pt}\ \par
\par\vs{-2pt}\ \par\cl{\bf\S1. Introduction}
\vs{4pt}\par
  A Lie algebra is called {\it complete} if its center is zero, and all its
derivations are inner. This definition of a complete Lie algebra
was given by Jacobson in 1962 [3]. In [12], Schenkman proved a
theorem, the so-called {\it derivation tower theorem}, that the
last term of the derivation tower of a centerless Lie algebra is
complete. It is also known that the holomorph of an abelian Lie
algebra is complete. Suggested by the derivation tower theorem,
complete Lie algebras occur naturally in the study of Lie algebras
and would provide interesting objects for investigation. In 1963,
Leger [6] presented some interesting examples of complete Lie
algebras. However since then, the study of complete Lie algebras
had become dormant for decades, partly due to the fact that the
structure theory of complete Lie algebras was not well developed.
In recent years, much progress has been obtained on the structure
theory of complete Lie algebras in finite dimensional case (see for example,
the references at the end of this paper).
\par It is known that a simple Lie algebra is not always complete
(e.g. in some infinite dimensional cases [14]). In this case, a
natural question is: {\it how big is the length of the derivation
tower of a simple Lie algebra?}
\par
In this paper, we answer the above question. Precisely, we prove
that the derivation algebra of a {\it centerless perfect} Lie
algebra (i.e., a Lie algebra $\fg$ with zero center and
$[\fg,\fg]=\fg$) of arbitrary dimension over any field of
arbitrary characteristic is complete, thus as a consequence, the
length of the derivation tower of any simple Lie algebra is $\le
1$.
\par
Let $\F$ be any field. First we recall the definitions of a
derivation and the holomorph of a Lie algebra $\fg$ over the field
$\F$. A {\it derivation} of a Lie algebra $\fg$ is an $\F$-linear
transformation $d:\fg\to\fg$ such that
$$d([x,y])=[d(x),y]+[x,d(y)]\;\;\;\;\for\;\;\;\;x,y\in\fg.
\eqno(1.1)
$$
We denote by $\derg$ the vector space of derivations of $\fg$,
which forms a Lie algebra with respect to the commutator of linear
transformations, called the {\it derivation algebra} of $\fg$.
Clearly, the space $\ad_{\fg}=\{\ad_x\,|\,x\in\fg\}$ of {\it inner
derivations} is an ideal of $\derg$. We call $\derg/\ad_{\fg}$ the
{\it outer derivation algebra} of $\fg$.
\par
 The {\it holomorph} ${\f h}(\fg)$ of a Lie algebra $\fg$ is the direct sum of
the vector spaces ${\f h}(\fg)=\fg\oplus\derg$ with the following
bracket
$$[(x,d),(y,e)]=([x,y]+d(y)-e(x),[d,e])\;\;\;\;\for\;\;\;\;x,y\in\fg,\,d,e\in\derg.\eqno(1.2)$$
An element $(x,d)$ of ${\f h}(\fg)$ is also written as $x+d$.
Obviously, $\fg$ is an ideal of ${\f h}(\fg)$ and ${\f
h}({\fg})/{\fg}\cong\derg$. Thus we write
$$
{\f h}(\fg)=\fg\semiprod\derg.\eqno(1.3)
$$
For a Lie algebra $\fg$, we denote by $C(\fg)$ the center of
$\fg$, i.e., $C(\fg)=\{x\in\fg\,|\,[x,\fg]=0\}$. \par The main
result of this paper is the following \vs{6pt}\par {\bf Theorem
1.1}. {\it Let $\f g$ be a perfect Lie algebra with zero center $($i.e.,
$[\fg,\fg]=\fg,\,C(\fg)=0)$.
Then we have
\par
(i) The derivation algebra ${\rm Der\,}{\f g}$ is complete.
\par
(ii) The holomorph ${\f h}({\f g})$ is complete if and only if the
center of outer derivation algebra is zero, i.e., $C({\rm
Der\,}{\f g}/ad_{\fg})=0$.} \vs{4pt}\par We shall prove Theorem
1.1 in Section 2, then we give some interesting example in Section
3.
\par\vs{-2pt}\ \par
\cl{\bf\S2. Proof of the main result} \vs{4pt}\par {\it Proof of
Theorem 1.1(i)}. Assume that $d\in C({\rm Der}{\f g})$. Then in
particular we have \linebreak[4]
$[d,\ad_x](y)=0$ for all $x,y\in {\f g}$. Thus
$d([x,y])=[x,d(y)]$. Hence by (1.1), $[d(x),y]=0$ for all
$x,y\in\fg$. Since $\fg$ has zero center, we obtain $d(x)=0$,
i.e., $d=0$.  Therefore
$$C({\rm Der\,}{\f g})=0.\eqno(2.1)$$
\par
Now we prove that all derivations of the Lie algebra $\derg$ are
inner. First we have \vs{4pt}\par {\bf Claim 1.} Let $ D\in {\rm
Der\,}(\derg)$. If $D(\ad_\fg)=0$, then $D=0$. \vs{2pt}\par Let
$d\in\derg,\,x\in\fg$. Note from (1.1) that in the Lie algebra
$\derg$, we have
$$
[d,\ad_x]=\ad_{d(x)}\in\ad_\fg.\eqno(2.2)$$ Using this, noting
that $D(d)\in\derg$ and the fact that $D(\ad_\fg)=0$, we have
$$
\begin{array}{ll}
\ad_{D(d)(x)}\!\!\!\!&=[D(d),\ad_x]\vs{4pt}\\
&=D([d,\ad_x])-[d,D(\ad_x)]\vs{4pt}\\
&=D(\ad_{d(x)})\vs{4pt}\\ &=0.\\ \end{array} \eqno(2.3)$$ Since
$\fg$ has zero center, (2.3) gives $D(d)(x)=0$ for all $x\in\fg$,
which means that $D(d)=0$ as a derivation of $\fg$. But $d$ is
arbitrary, we obtain $D=0$. This proves the claim. \vs{4pt}\par
Since $\fg$ is perfect, for any $x\in\fg$, we can write $x$ as
$$x=\sum_{i\in I}[x_i,y_i]\;\;\;\;\mbox{for \
some}\;\;\;\;x_i,y_i\in\fg,\eqno(2.4)$$ and for some finite index
set $I$. Then
$$
\ad_x=\sum_{i\in I}[\ad_{x_i},\ad_{y_i}].\eqno(2.5)
$$
Then for any $D\in{\rm Der}(\derg)$, we have
$$
\begin{array}{ll}D(\ad_x)\!\!\!\!&=\dis\sum_{i\in I}D([\ad_{x_i},\ad_{y_i}])\vs{4pt}\\
&\dis=\sum_{i\in I} ([D(\ad_{x_i}),\ad_{y_i}]+[\ad_{x_i},
D(\ad_{y_i})]).\\ \end{array}\eqno(2.6)$$ Let
$d_i=D(\ad_{x_i}),\,e_i=D(\ad_{y_i})\in\derg$. Then by (2.2), we
have
$$\begin{array}{ll}D(\ad_x)\!\!\!\!&\dis=\sum_{i\in I}(\ad_{d_i(y_i)}-\ad_{e_i(x_i)})\vs{4pt}\\ &
=\ad_{_{\sc\sum_{i\in I}(d_i(y_i)-e_i(x_i))}}\in \ad_\fg.\\
\end{array}\eqno(2.7)$$ This means that $D(\ad_x)=\ad_y$ for some
$y\in\fg$. Since $C({\fg})=0$, such $y$ is unique. Thus
$d:x\mapsto y$ defines a linear transformation of $\fg$ such that
$$
D(\ad_x)=\ad_{d(x)}.\eqno(2.8)$$ For $x,y\in\fg$, we have
$$\begin{array}{ll}
\ad_{d([x,y])}\!\!\!\!&=D(\ad_{[x,y]})\vs{4pt}\\
&=D([\ad_x,\ad_y])\vs{4pt}\\
&=[D(\ad_x),\ad_y]+[\ad_x,D(\ad_y)]\vs{4pt}\\
&=[\ad_{d(x)},\ad_y]+[\ad_x,\ad_{d(y)}]\vs{4pt}\\
&=\ad_{[d(x),y]+[x,d(y)]}.\\ \end{array}\eqno(2.9)
$$
This and the fact that $C(\fg)=0$ mean that
$d([x,y])=[d(x),y]+[x,d(y)]$, i.e., $d\in\derg$. Then (2.8) gives
that $D(\ad_x)=\ad_{d(x)}=[d,\ad_x]$ for all $x\in\fg$, i.e.,
$$(D-\ad_d)(\ad_\fg)=0.\eqno(2.10)$$
By Claim I, we have $D-\ad_d=0$, i.e., $D=\ad_d$ is an inner
derivation on $\derg$. This together with (2.1) proves Theorem
1.1(i). \vs{6pt}
\par {\it Proof of Theorem 1.1(ii)}. ``$\Longleftarrow$'': First we prove the sufficiency.
So suppose $C(\derg/\ad_\fg)=0$. We want to prove that ${\f h}(\fg)$ is
complete.
\par First we prove $C({\f h}(\fg))=0$. Suppose $h=x+d\in C({\f
h}(\fg))$ for some $x\in\fg,\,d\in\derg$. Letting $y=0$ in (1.2),
we obtain $[d,e]=0$ for all $e\in\derg$, i.e., $d\in C(\derg)=0$.
Then $h=x\in C({\f h}(\fg))\cap\fg\subset C(\fg)=0$, i.e., $h=0$.
Thus $C({\f h}(\fg))=0$. Then as in the proof of (2.1), we have
$$C({\rm Der}({\f h}(\fg)))=0.\eqno(2.11)$$ Now let $\DD\in {\rm
Der}({\f h}({\f g}))$. For any $x\in\fg\subset{\f h}(\fg)$ written
in the form (2.4), since $\fg$ is an ideal of ${\f h}(\fg)$, we
have
$$\DD(x)=\sum_{i\in
I}([\DD(x_i),y_i]+[x_i,\DD(y_i)])\in\fg.\eqno(2.12)$$ Thus
$d'=\DD|_{\fg}:\fg\to\fg$ is a derivation of $\fg$, i.e.,
$d'\in\derg\subset{\f h}(\fg)$. Let $\DD_1=\DD-\ad_{d'}\in{\rm
Der}({\f h}(\fg))$. Then $\DD_1(x)=\DD(x)-[d',x]=\DD(x)-d'(x)=0$,
i.e.,
$$\DD_1|_{\fg}=0.\eqno(2.13)$$
\par
For any $d\in \derg\subset{\f h}(\fg)$, by (1.3), we can write
$\DD_1(d)\in{\f h}(\fg)$ as
$$\DD_1(d)=x_d+d_1\;\;\;\;\mbox{ for \ some}\;\;\;\;x_d\in\fg,\,d_1\in\derg.\eqno(2.14)$$
Then in ${\f h}(\fg)$, by (1.2), for any $y\in\fg$, we have
$$\begin{array}{ll}
[x_d,y]+d_1(y)\!\!\!\!&=[x_d+d_1,y]\vs{4pt}\\
&=[\DD_1(d),x_d]\vs{4pt}\\
&=\DD_1([d,x_d])-[d,\DD_1(x_d)]\vs{4pt}\\
&= \DD_1(d(x_d))-[d,\DD_1(x_d)]\vs{4pt}\\
&=0,
\\ \end{array}
\eqno(2.15)$$%
 where the last equality follows from (2.13). This
means that $d_1=-\ad_{x_{_{\ssc d}}}$ and so
$\DD_1(d)=x_d-\ad_{x_{_{\ssc d}}}$. For $d,e\in\derg$, we have $$
\begin{array}{ll}x_{[d,e]}-\ad_{x_{_{\ssc[d,e]}}}\!\!\!\!&=
\DD_1([d,e])\vs{4pt}\\ &=[\DD_1(d),e]+[d,\DD_1(e)] \vs{4pt}\\
&=(-e(x_d)+d(x_{e}))+([-\ad_{x_{_{\ssc d}}},e]+[d,-\ad_{x_{_{\ssc
e}}}]),
\\ \end{array}\eqno(2.16)$$
where the last equality follows from (1.2) and the fact that
$\DD_1(d)=x_d-\ad_{x_{_{\ssc d}}}$. So
$$
\ad_{x_{_{\ssc[d,e]}}}=[\ad_{x_{_{\ssc d}}},e]+[d,\ad_{x_{_{\ssc
e}}}], \eqno(2.17)$$ i.e., the linear transformation
$D:\derg\to\derg$ defined by
$$
D(d)=\ad_{x_{_{\ssc d}}}\;\;\;\;\for\;\;\;\;d\in\derg,
\eqno(2.18)$$ is a derivation of $\derg$. Since $\derg$ is
complete, there exists $d''\in\derg$ such that
$$D=\ad_{d''}.\eqno(2.19)$$
For any $d\in\derg$, we have
$$[d'',d]=D(d)=\ad_{x_{_{\ssc d}}}\in\ad_\fg,
\eqno(2.20)$$ i.e., $d''+\ad_\fg\in C(\derg/\ad_\fg)=0$. Thus
$d''\in\ad_\fg$. Therefore there exists $y\in\fg$ such that
$$
d''=\ad_y.\eqno(2.21)$$ Note that $y-\ad_y\in{\f h}(\fg)$. Let
$\DD_2=\DD_1-\ad_{y-\ad_y}$. Then for any $x\in\fg\subset{\f
h}(\fg)$, we have
$$\begin{array}{ll}\DD_2(x)\!\!\!\!&=\DD_1(x)-[y-\ad_y,x]
\vs{4pt}\\ &=-[y,x]+\ad_y(x)\vs{4pt}\\
&=-[y,x]+[y,x]\vs{4pt}\\ &=0,\\ \end{array}\eqno(2.22)$$ where the
second equality follows from (2.13) and (1.2), and for any
$d\in\derg\subset{\f h}(\fg)$,
$$
\begin{array}{ll}
\DD_2(d)\!\!\!\!&=\DD_1(d)-[y-\ad_y,d]\vs{4pt}\\ &
=(x_d-\ad_{x_{_{\ssc d}}})-(-d(y)-[\ad_y,d])\vs{4pt}\\
&=(x_d+d(y))-(\ad_{x_{_{\ssc d}}}-[\ad_y,d])\vs{4pt}\\
&=0,\end{array}\eqno(2.23)$$ because
$-\ad_{d(y)}=[\ad_y,d]=[d'',d]=\ad_{x_{_{\ssc d}}}$ by (2.20) and
(2.21), and $x_d=-d(y)$ (since $\fg$ has zero center). Thus (2.22)
and (2.23) show that $\DD_2=0$ and so $\DD$ is inner. This and
(2.11) prove that ${\f h}(\fg)$ is complete. \vs{6pt}\par
$``\Longrightarrow$'': Now we prove the necessity. So assume that
${\f h}(\fg)$ is complete. Suppose conversely
$C(\derg/\ad_\fg)\ne0$. Then there exists $D\in\derg$ such that
$$D\notin\ad_\fg\;\;\;\;\mbox{but}\;\;\;\;[D,\derg]\subset\ad_\fg.\eqno(2..24)$$
Then for any $d\in\derg$, there exists $x_d\in\fg$ such that
$[D,d]=\ad_{x_{_{\ssc d}}}$. Such $x_d$ is unique since
$C(\fg)=0$. Using the fact that $[D,[d,e]]=[[D,d],e]+[d,[D,e]]$,
we obtain
$$
x_{[d,e]}=-e(x_d)+x_d(e)\;\;\;\;\for\;\;\;\;d,e\in\derg.
\eqno(2.25)$$ We define a linear transformation $\DD:{\f
h}(\fg)\to{\f h}(\fg)$ as follows: $\DD|_\fg=0$, and for
$d\in\derg\subset{\f h}(\fg)$, we define
$$
\DD(d)=x_d-\ad_{x_{_{\ssc d}}}\in \fg\semiprod\derg={\f h}(\fg).
\eqno(2.26)$$ From this definition, we have
$$[\DD(d),x]=0\;\;\;\;\for\;\;\;\;x\in\fg\subset{\f
h}(\fg),\,d\in\derg\subset{\f h}(\fg).\eqno(2.27)$$ Then for any
$h=x+d,\,h'=y+e\in{\f h}(\fg)$, by (1.2) and the fact that
$\DD_\fg=0$, we have
$$
\begin{array}{ll}
\DD([h,h'])\!\!\!\!&=\DD([d,e])\vs{4pt}\\
&=x_{[d,e]}-\ad_{x_{_{\ssc[d,e]}}} \vs{4pt}\\ &=
(-e(x_d)+d(x_e))-([\ad_{x_{_{\ssc d}}},e]+[d,\ad_{x_{_{\ssc e}}}]) \vs{4pt}\\
&=(-e(x_d)-[\ad_{x_{_{\ssc d}}},e])+(d(x_e)+[d,\ad_{x_{_{\ssc
e}}}])\vs{4pt}\\ &=
[\DD(d),e]+[d,\DD(e)]\vs{4pt}\\ &=
[\DD(h),h']+[h,\DD(h')].\\
\end{array}
\eqno(2.28)$$ Thus $\DD$ is a derivation of ${\f h}(\fg)$. Since
${\f h}(\fg)$ is complete, $\DD$ is inner, therefore, there exists
$h=y+e\in{\f h}(\fg)$ such that $\DD=\ad_h$. For any $x\in\fg$,
since $\DD|_\fg=0$, we have
$$
\begin{array}{ll}
0\!\!\!\!&=\DD(x)\vs{4pt}\\ &=[h,x]\vs{4pt}\\ &=[y,x]+e(x),\\
\end{array}
\eqno(2.29)$$ i.e., $e=-\ad_y$. Then by (2.26),
$$
\begin{array}{ll}
x_d-\ad_{x_{_{\ssc d}}}\!\!\!\!&=\DD(d)\vs{4pt}\\
&=[h,d]\vs{4pt}\\ &=-d(y)-[\ad_y,d]\vs{4pt}\\
&=-d(y)+\ad_{d(y)}.\\ \end{array}\eqno(2.30)$$ Hence
$\ad_{x_{_{\ssc d}}}=-\ad_{d(y)}$. Then for any $d\in\derg$,
$$
\begin{array}{ll}
[D,d]\!\!\!\!&=\ad_{x_{_{\ssc d}}}\vs{4pt}\\
&=-\ad_{d(y)}\vs{4pt}\\
&=[\ad_y,d].\\ \end{array} \eqno(2.31)$$ Since $\derg$ has zero
center, (2.31) implies that $D=\ad_y\in\ad_\fg$, a contradiction
with (2.24). Thus $C(\derg/\ad_\fg)=0$, and the proof of Theorem
1.1(ii) is complete.
\par\vs{-2pt}\ \par
\cl{\bf\S3. Some  Examples} \vs{4pt}\par {\bf Example 3.1}. An
$n\times n$ matrix $q=(q_{ij})$ over a field $\F$ of
characteristic $0$ such that
$$q_{ii}=1 \ \ \mbox{and} \ \ q_{ji}=q_{ij}^{-1},\eqno(3.1)$$
is called a {\it quantum matrix}. The {\it quantum torus}
$$\F_q=\F_q[t_1^{\pm},\cdots,t_n^{\pm}],\eqno(3.2)$$
determined by a quantum matrix $q$ is defined as the associative
algebra over $\F$ with $2n$ generators
$t_1^{\pm},\cdots,t_n^{\pm}$, and relations
$$t_it_i^{-1}=t_i^{-1}t_i=1 \ \ \mbox{and} \ \ t_jt_i=q_{ij}t_it_j,\eqno(3.3)$$ for all $1\leq i,j\leq n$.
\par
For any $a=(a_1,\cdots,a_n)\in \Z^n$, we denote
$t^a=t_1^{a_1}\cdots t_n^{a_n}$. For any $a,b\in \Z^n$, we define
$$\sigma(a,b)=\prod_{1\leq i,j\leq n}q_{j,i}^{a_jb_i} \ \ \mbox{and} \ \
f(a,b)=\prod_{i,j=1}^nq_{j,i}^{a_jb_i}.\eqno(3.4)$$ Then we have
$$t^at^b=\sigma(a,b)t^{a+b},\ \ t^at^b=f(a,b)t^bt^a\mbox{ \ and
\ }f(a,b)=\sigma(a,b)\sigma(b,a)^{-1}.%
\eqno(3.5)$$%
Note that the commutator of monomials in $\F_q$ satisfies
$$[t^a,t^b]=(\sigma(a,b)-\sigma(b,a))t^{a+b}=\sigma(b,a)(f(a,b)-1)t^{a+b},\eqno(3.6)$$
for all $a,b\in \Z^n$. Define the {\it radical} of $f$, denoted by
${\rm rad}(f)$, by
$${\rm rad}(f)=\{a\in \Z^n\,|\,f(a,b)=1 \ \ \mbox{for \ all}\ \ b\in \Z^n\}.\eqno(3.7)$$
It is clear from (3.4) that ${\rm rad}(f)$ is a subgroup of
$\Z^n$. The center of $\F_q$ is $C(\F_q)=\sum_{a\in {\rm
rad}(f)}\F t^a$ and the Lie algebra $\F_q$ has the Lie ideals
decomposition
$$\F_q=C(\F_q)\oplus [\F_q,\F_q].\eqno(3.8)$$
By (3.8), we obtain that the Lie algebra $[\F_q,\F_q]\cong
\F_q/C(\F_q)$ is perfect and has zero center. From Theorem 1.1, we
obtain the following result. \vs{6pt}\par {\bf Corollary 3.2}.
{\it Let $\F_q$ be a quantum torus as above. Then the derivation
algebra\linebreak[4] ${\rm Der}([\F_q,\F_q])$ of $[\F_q,\F_q]$ is complete.}
\vs{6pt}\par {\bf Example 3.3}. A Lie algebra $\fg$ is called a
{\it symmetric self-dual Lie algebra} if $\fg$ is endowed with a
nondegenerate invariant symmetric bilinear form $B$. For any
subspace $V$ of $\fg$, we define $V^{\bot}=\{x\in
{\fg}\,|\,B(x,y)=0,\forall\, y \in V\}$. Then we can easily check
that $[{\frak g},{\frak g}]^{\bot}=C(\fg)$, where $C(\fg)$ is the
center of $\fg$. Then we have \vs{6pt}\par {\bf Corollary 3.4}.
{\it Assume that $\fg$ is a symmetric self-dual Lie algebra with
zero center. Then the Lie algebra $\derg$ is complete.} \vs{6pt}
\par
{\bf Remark 3.5}. By [18], we know that the structure of a perfect
symmetric self-dual Lie algebra is not clear even in the finite
dimension case.
\par\vs{-2pt}\ \par
\cl{\bf References}\vs{4pt}\par \ni\hi3.5ex\ha1[1] S.~Benayadi,
``The perfect Lie algebras without center and outer derivations,''
{\it Annals de la Faculte des Sciences de Toulouse} {\bf 5}
(1996), 203-231.
\par \ni\hi3.5ex\ha1[2] R.~Carles, ``Construction des algebres de Lie completes,''
{\it C.~R.~Acad.~Sci.~Paris Ser.~I.~Math.} {\bf 318} (1994),
171-174.
\par \ni\hi3.5ex\ha1[3] N.~Jacobson, {\it Lie algebras}, Wiley (Interscience), New York, 1962.
\par \ni\hi3.5ex\ha1[4] C.~P.~Jiang, D.~J.~Meng and S.~Q.~Zhang, ``Some complete Lie
algebras,'' {\it J. Algebras} {\bf 186} (1996), 807-817.
\par \ni\hi3.5ex\ha1[5] V.~G.~Kac, {\it Infinite dimensional Lie algebras}, 3rd ed, Cambridge
Univ.~Press, 1990.
\par \ni\hi3.5ex\ha1[6] G.~Leger, ``Derivations of Lie algebras III,'' {\it Duke
Math.~J.}~{\bf 30} (1963), 637--645.
\par \ni\hi3.5ex\ha1[7] D.~J.~Meng, ``Some results on complete Lie
algebras,'' {\it Comm.~Algebra} {\bf 22} (1994), 5457-5507.
\par \ni\hi3.5ex\ha1[8] D.~J.~Meng, ``Complete Lie algebras and Heisenberg
algebras,'' {\it Comm.~Algebra} {\bf 22} (1994), 5509-5524.
\par \ni\hi3.5ex\ha1[9] D.~J.~Meng, S.~P.~Wang, ``On the construction of complete Lie
algebras,'' {\it J.~Algebra} {\bf 176} (1995), 621-637.
\par \ni\hi4ex\ha1[10] D.~J.~Meng, L.~S.~Zhu, ``Solvable complete Lie algebras I,'' {\it Comm.~Algebra}
{\bf 24} (1996), 4181-4197.
\par \ni\hi4ex\ha1[11] E.~L.~Stitzinger, ``On Lie algebras with only inner
derivations,'' {\it J.~Algebra} {\bf 105} (1987), 341--343.
\par \ni\hi4ex\ha1[12] E.~V.~Schenkman, ``A theory of subinvariant Lie algebras,'' {\it Amer.~J.~Math.}
{\bf 73} (1951), 453--474.
\par \ni\hi4ex\ha1[13] S.~T\^og\^o, ``Outer derivations of Lie algebras,'' {\it Trans.~Amer.~Math.~Soc.}
{\bf 128} (1967), 264--276.
\par \ni\hi4ex\ha1[14] Q.~Yang, D.~J.~Meng, ``Derivation algebras of
simple Lie algebras,'' {\it Dongbei Shuxue}, In press.
\par \ni\hi4ex\ha1[15] L.~S.~Zhu, D.~J.~Meng, ``Derivations of generalized Kac-Moody
algebras,'' {\it Chinese Annal of Mathematics} {\bf 17A:3} (1996),
271-176 (in chinese).
\par \ni\hi4ex\ha1[16] L.~S.~Zhu, D.~J.~Meng, ``The complete Lie algebras with maximal
rank nilpotent radicals,'' {\it Chinese Science Bulletin} {\bf 44}
(1999), 4.
\par \ni\hi4ex\ha1[17] L.~S.~Zhu, D.~J.~Meng, ``Complete Lie algebras II,'' {\it Alg.~Colloq.}
{\bf 5} (1998), 289-296.
\par \ni\hi4ex\ha1[18] L.~S.~Zhu, D.~J.~Meng, ``A class of Lie algebras with a symmetric invariant non-degenerate
bilinear form,'' {\it Acta Math. Sinica} {\bf 43} (2000),
1119-1126.

\end{document}